\documentclass{article}
\usepackage{amsmath}
\usepackage{amsfonts,amsthm}
\usepackage{amssymb,enumerate,enumitem,verbatim}
\usepackage[pdftex]{graphicx}

\newtheorem{lemma}{Lemma}[section]
\newtheorem{corollary}[lemma]{Corollary}
\newtheorem{theorem}[lemma]{Theorem}

\theoremstyle{definition}
\newtheorem*{conjecture}{Conjecture}
\newtheorem*{acknowledgements}{Acknowledgements}
\newtheorem*{defn}{Definition}

\theoremstyle{remark}
\newtheorem*{rmk}{Remark}

\newtheorem*{example}{Example}
\newtheorem*{nte}{Note}

\global\long\def\a{\alpha}

\global\long\def\e{\varepsilon}

\global\long\def\N{\mathbb{N}}

\global\long\def\GG{\mathcal{G}}

\global\long\def\TT{\mathcal{T}}

\global\long\def\PP{\mathcal{P}}

\global\long\def\re{\begin{rmk}}
\global\long\def\mark{\end{rmk}}
\global\long\def\ex{\begin{example}}
\global\long\def\ple{\end{example}}
\global\long\def\no{\begin{nte}}
\global\long\def\ted{\end{nte}}
\global\long\def\en{\begin{compactenum}}
\global\long\def\um{\end{compactenum}}
\global\long\def\li{\begin{compactitem}}
\global\long\def\st{\end{compactitem}}
\global\long\def\de{\begin{defn}}
\global\long\def\fn{\end{defn}}
\global\long\def\cor{\begin{corollary}}
\global\long\def\ary{\end{corollary}}
\global\long\def\lem{\begin{lemma}}
\global\long\def\ma{\end{lemma}}
\global\long\def\arr{\begin{array}}
\global\long\def\ay{\end{array}}
\global\long\def\pr{\begin{proof}}
\global\long\def\oof{\end{proof}}

\newif\ifdraft
\draftfalse 

\newif\ifdone
\donefalse

\title{Embedding bounded degree spanning trees in random graphs}

\author{Richard Montgomery\footnote{Department of Pure Mathematics and Mathematical Statistics, Centre for Mathematical Sciences, Wilberforce Road, Cambridge, CB3 0WB, UK. r.h.montgomery@dpmms.cam.ac.uk}}

\begin{document}

\maketitle
\begin{abstract} 
We prove that if a tree $T$ has $n$ vertices and maximum degree at most $\Delta$, then a copy of $T$ can almost surely be found in the random graph $\mathcal{G}(n,\Delta\log^5 n/n)$.
\end{abstract}

\section{Introduction}\label{intro}

The threshold for a binomial random graph $\GG(n,p)$ to contain a spanning tree is well known, coinciding trivially with the threshold for the graph to be connected at $p=\log n/ n$. Which spanning trees may we expect to appear at, or close to, this threshold? Around twenty years ago, Kahn made the following natural conjecture.

\begin{conjecture}
For every fixed $\Delta>0$, there is some constant $C$ such that if $T$ is a tree on $n$ vertices with maximum degree at most $\Delta$, then the random graph $\GG(n,C\log n/ n)$ almost surely\footnote{More formally, given any sequence of such trees $T_n$, with $|T_n|=n$, $\PP(T_n\subset \GG(n,C\log n/ n))\to 1$.} contains a copy of $T$.
\end{conjecture}

Alon, Krivelevich and Sudakov~\cite{AKS07} studied the related problem of embedding almost spanning bounded degree trees in a random graph. They showed that, for sufficiently large $c=c(\e,\Delta)$, the random graph $\GG(n,c/n)$ almost surely contains every tree with maximum degree at most $\Delta$ consisting of at most $(1-\e)n$ vertices. This is best possible up to the constant $c$. An improved bound for this constant was proved by Balogh, Csaba, Pei and Samotij~\cite{BCPS10}, and a resilience version of the result was obtained by Balogh, Csaba and Samotij~\cite{BCS11}.

Where a tree $T$ has $n$ vertices and maximum degree $\Delta$, Krivelevich \cite{MK10} showed that, for any $\e>0$, the random graph $\GG(n,\max\{(40/\e)\Delta\log n,n^\e\}/n)$ almost surely contains a copy of $T$.  When $\Delta=n^{\Theta(1)}$ this result is essentially tight, but for smaller $\Delta$ the situation is less clear. For example, where $\Delta$ is constant this leaves a large gap to the conjecture stated above. The theorem below reduces this gap.

\begin{theorem}\label{mainexpo1} 
If $T$ has $n$ vertices and maximum degree $\Delta$, then the random graph $G=\GG(n,p)$ with $p=\Delta\log^5 n/ n$ almost surely contains a copy of $T$.
\end{theorem}

The techniques used here to prove Theorem \ref{mainexpo1} can be applied more carefully to further reduce the probability $p$ required. Even with care, however, and some development of the techniques, the method reaches a natural barrier around $p=\Delta\log^2 n/ n$. Additional ideas, resulting in more technicalities, are required to break through this barrier and prove the conjecture stated earlier. In anticipation of a paper by the author achieving this, we do not attempt to minimize the probability in the method presented here (see Section \ref{future}). Instead, we concentrate on giving a clear proof of Theorem \ref{mainexpo1}, which already improves distinctly on what is known in the case of small $\Delta$. We also record a path covering result that will prove useful elsewhere~\cite{selflove2}.

In~\cite{MK10}, Krivelevich considered trees differently according to whether they contained many leaves or many paths without any branch vertices. Such a path in a tree, where all the interior vertices have degree 2 in the parent tree, we call a \emph{bare path}. Spanning trees with many leaves may be embedded using the following argument outlined by Alon, Krivelevich and Sudakov~\cite{AKS07}. If a tree $T$ contains $\e n$ leaves, for fixed $\e>0$, and has maximum degree $\Delta$, we consider the subtree $S$ gained from deleting the leaves of $T$. The tree $S$ can be embedded in a random graph using the almost spanning tree result stated earlier. To complete the embedding of $T$ we must match the vertices in the copy of $S$ which need leaves added to the correct number of remaining vertices in $G$. By revealing more edges, such a matching can be found to complete the embedding. This scheme works almost surely when $pn\geq (40/\e)\Delta\log n/ n$~\cite{MK10}.

As Krivelevich showed~\cite{MK10}, trees with fewer than $n/4k$ leaves must have, by way of compensation, $n/4k$ disjoint bare paths of length $k$ (Lemma \ref{split}). The forest gained by removing these paths may be embedded in a random graph using the almost spanning tree result. The remaining vertices then need to be covered by disjoint paths of length $k$, connecting specific pairs of vertices to complete the desired spanning tree. Taking $k$ to be constant, Krivelevich used a modification of the very strong results by Johansson, Kahn and Vu~\cite{JKV08} to complete the embedding, revealing more edges in the random graph to find the missing paths.

We follow this outline, but in the final step cover the remaining vertices by paths using a different method. This allows longer paths to be used, and the covering to be done at a much lower probability. We use a version of \emph{absorption}, given as a general method by R\"odl, Ruci\'nski and Szemer\'edi~\cite{RRS08}. In their terminology, we use a natural construction for \emph{absorbers}, similar in spirit to that used elsewhere by K\"uhn and Osthus~\cite{KO12}, and Allen, B\"{o}ttcher, Kohayakawa and Person~\cite{ABKP13}. Our methods differ from previous implementations of this method in two main ways. Using methods in Section \ref{paths} to create paths between vertices, we create our absorbers in a very sparse random graph. We also construct our \emph{reservoir} using a system of paths, rather than a single long path.

In Section \ref{tools} we recall some previous results, and construct graphs with a certain useful property. In Section \ref{paths} we give a construction for paths in expander graphs. In Section \ref{main} we prove Theorem \ref{mainexpo1}, before recording a useful variation of our main path covering theorem in Section \ref{technical}. In Section \ref{future} we discuss some further work.

We use $\log n$ for the natural logarithm, and omit rounding symbols where they are not crucial.

\section{Preliminaries}\label{tools}

\subsection{Notation}
For a graph $G$, $V(G)$ is the vertex set of $G$, and $|G|=|V(G)|$. For a subset $U\subset V(G)$, let $G[U]$ be the subgraph of $G$ induced on the set $U$. Given two disjoint sets $X$ and $Y$, $e_G(X,Y)$ is the number of edges between $X$ and $Y$. For a subset $U\subset V(G)$ and vertices $x,y\in U$, $d_U(x,y)$ is the length of a shortest $x,y$-path in $G[U]$, if such a path exists. The set of neighbours of a vertex $v$ is denoted by $N(v)$, and the neighbourhood of a vertex $v$ in the set $A\subset V(G)$ is $N(v,A)=N(v)\cap A$. The neighbourhood of a vertex set $S\subset V(G)$ is $N(S)=(\cup_{v\in S}N(v))\setminus S$ and $N(S,A)=N(S)\cap A$. Where multiple graphs are involved we use $N_G(v,A)$ to emphasise the graph currently considered. We say a path with $l$ vertices has \emph{length} $l-1$, and call a path $P$ an $x,y$-path if the vertices $x$ and $y$ have degree 1 in $P$.

Given two disjoint vertex sets $A$ and $B$, a \emph{$d$-matching from $A$ into $B$} is a collection of disjoint sets $\{X_a\subset N(a,B):a\in A\}$ so that, for each $a\in A$, $|X_a|=d$. As is well known, such matchings can be found by showing that Hall's generalised matching condition holds. For details on this, and other standard notation, see Bollob\'as~\cite{bollo1}.

\subsection{Graph Expansion and Trees}
We will establish graph expansion properties in the random graph, before using the properties to embed the trees.

\de Let $n\in \mathbb{N}$ and $d>0$. A graph $G$ is an \emph{$(n,d)$-expander} if $|G|=n$ and $G$ satisfies the following two conditions.
\begin{enumerate}
\item $|N_G(X)|\geq d|X|$ for all sets $X\subset V(G)$ with $1\leq |X|<\lceil \frac{n}{2d}\rceil$.
\item $e_G(X,Y)>0$ for all disjoint $X,Y\subset V(G)$ with $|X|=|Y|=\lceil\frac{n}{2d}\rceil$.
\end{enumerate}
\fn

\de
Let $\TT(n,\Delta)$ be the class of all trees on $n$ vertices with maximum degree at most $\Delta$.
\fn

As shown by Balogh, Csaba, Pei and Samotij~\cite{BCPS10}, almost spanning trees can be found in expander graphs using a theorem by Haxell~\cite{PH01}. We will use the following formulation due to Johannsen, Krivelevich and Samotij~\cite{JKS12}.
\begin{theorem}\label{almost}
Let $n,\Delta\in \mathbb{N}$, let $d\in \mathbb{R}^+$ with $d\geq 2\Delta$ and let $G$ be a $(n,d)$-expander. Given any $T\in \TT(n-4\Delta\lceil\frac{n}{2d}\rceil,\Delta)$, we can find a copy of $T$ in the graph $G$.
\end{theorem}

For Section \ref{paths} we will also find the following formulation by Balogh, Csaba, Pei and Samotij useful~\cite{BCPS10}.
\begin{theorem}\label{almost2}
Let $\Delta,m,M\in \N$. If $H$ is a non-empty graph such that for every $X\subset V(H)$, if $0<|X|\leq m$, then $N_H(X)|\geq \Delta|X|+1$ and, if $|X|=m$, then $|N_H(X)|\geq 2\Delta m+M$, then $H$ contains every tree in $\TT(M,\Delta)$.
\end{theorem}

For our constructions we wish to partition our expander graph so that all vertex subsets expand well into each vertex subset in the partition. The following partition lemma by Johanssen, Krivelevich and Samotij~\cite{JKS12} permits this. For convenience, we state a slightly stronger form of their lemma, using the following definition, but the proof follows identically.

\de
For a graph $G$ and a set $W\subset V(G)$, we say $G$ \emph{$d$-expands} into $W$ if
\begin{enumerate}
\item $|N_G(X, W)|\geq d|X|$ for all $X\subset V(G)$ with $1\leq |X|<\lceil \frac{|W|}{2d}\rceil$, and,
\item $e_G(X,Y)>0$ for all disjoint $X,Y\subset V(G)$ with $|X|=|Y|=\lceil\frac{|W|}{2d}\rceil$.
\end{enumerate}
\fn

\lem[\cite{JKS12}]\label{splitexpand} There exists an absolute constant $n_0\in \N$ such that the following statement holds. Let $k,n\in \N$ and $d\in \mathbb{R}^+$ satisfy $n\geq n_0$ and $k\leq \log n$. Furthermore, let $m,m_1,\ldots,m_k\in \N$ satisfy $m=m_1+\ldots+m_k$ and let $d_i:=\frac{m_i}{5m}d$ satisfy $d_i\geq 2\log n$ for all $i\in \{1,\ldots,k\}$. Then, for any graph $G$ which $d$-expands into a vertex set $W$, with $|W|=m$, the set $W$ can be partitioned into $k$ disjoint sets $W_1,\ldots, W_k$ of sizes $m_1,\ldots,m_k$ respectively, such that, for each $i$, $G$ $d_i$-expands into $W_i$.
\ma

In order to find expansion properties in the binomial random graph we will use the following lemma by Johannsen, Krivelevich and Samotij~\cite{JKS12}.

\lem\label{randomexpand} Let $d:\N\to\mathbb{R}^+$ satisfy $d\geq 3$. Then $\GG(n,7dn^{-1}\log n)$ is almost surely an $(n,d)$-expander.
\ma

As mentioned in Section \ref{intro}, we will use the following lemma by Krivelevich~\cite{MK10} to divide $\TT(n,\Delta)$ into trees with many leaves and trees with many disjoint bare paths.

\lem\label{Kleaves} Let $k,l,n>0$ be integers. Let $T$ be a tree on $n$ vertices with at most $l$ leaves. Then $T$ contains a collection of at least $n/(k+1)-(2l-2)$ vertex disjoint bare paths of length $k$ each.
\ma
We will employ this lemma in the following form, letting $l=n/4k$.

\begin{corollary}\label{split} For any integers $n,k>2$, a tree on $n$ vertices either has at least $n/4k$ leaves or a collection of at least $n/4k$ vertex disjoint bare paths of length $k$ each.
\end{corollary}

We will also use the following lemma by Krivelevich~\cite{MK10} to almost surely attach leaves onto a partial embedding of a tree.
\lem\label{matchymatchy} Let $0<d_1,\ldots,d_k$ be integers satisfying: $d_i\leq \Delta$, $\sum^k_{i=1}d_i=l$. Let $A=\{a_1,\ldots,a_k\}$, $B$ be disjoint sets of vertices with $|B|=l$. Let $G$ be a random bipartite graph with sides $A$ and $B$, where each pair $(a,b)$, $a\in A$, $b\in B$, is an edge of $G$ with probability $p$, independent of other pairs. If $p\geq 2\Delta\log l/ l$, then, almost surely as $l\to\infty$, the random graph $G$ contains a collection $S_1,\ldots, S_k$ of vertex disjoint stars such that $S_i$ is centered at $a_i$ and has the remaining $d_i$ vertices in $B$.
\ma

\subsection{A random graph construction}
For Theorem \ref{mainexpo1}, we require a bipartite graph with a certain `resilient matching' condition. Here we verify that such graphs exist.
\lem \label{flexiblematching}
There is a constant $n_0$, such that for every $n\geq n_0$ with $3|n$, there exists a bipartite graph $H$ on vertex classes $X$ and $Y\cup Z$ with $|X|=n$, $|Y|=|Z|=2n/3$, and maximum degree $40$, so that the following is true. Given any subset $Z'\subset Z$ with $|Z'|=n/3$, there is a matching between $X$ and $Y\cup Z'$.
\ma
\pr
Let $m=2n/3$, and take two vertex sets $X_1$, $Y$ with $|X_1|=|Y|=2n/3$. Independently, place 20 random matchings between $X_1$ and $Y$ and let the graph $G$ be the union of these matchings. Given two sets $A\subset X_1$ and $B\subset Y$, and a random matching $M$, the probability that $N_M(A)\subset B$ is $\binom{|B|}{|A|}\binom{m}{|A|}^{-1}$. Thus, for each $t\leq m/4$, the probability there is some set $A\subset X_1$, $|A|=t$ with $|N_G(A)|<2t$, is at most
\[
\binom{m}{t}\binom{m}{2t}\left(\binom{2t}{t}\binom{m}{t}^{-1}\right)^{20}
\leq \left(\frac{e^3}{4}\left(\frac mt\right)^{3}\left(\frac{2t}{m}\right)^{20}\right)^t
= \left(\frac{8e^3}{4}\left(\frac{2t}{m}\right)^{17}\right)^t.
\]
Looking separately at the cases $t\leq \log m$, and $\log m\leq t\leq m/4$, we see that the sum of these probabilities tends to 0 as $m\to\infty$. An identical calculation shows that, almost surely, $|N(B)|\geq 2|B|$ for every $B\subset Y$ with $|B|\leq m/4$.
The probability there are no two sets $A\subset X_1$ and $B\subset Y$, with $|A|=|B|=\lceil m/4\rceil$ and $e_G(A,B)=0$ is at most
\[
\binom{m}{m/4}^2\left(\binom{3m/4}{m/4}\binom{m}{m/4}^{-1}\right)^{20}\leq (4e)^{m/2}\left(\frac{3}{4}\right)^{20m/4}\leq \left(\frac12\right)^{m/4}. 
\]
Therefore, almost surely, all three of these properties mentioned hold. Let $n$, and hence $m$, be sufficiently large that some graph with these properties must exist, and fix such a graph $G$. Duplicate the vertices in $Y$ to get $Z$, and then duplicate $m/2$ of the vertices in $X_1$ to get the set $X_2$. Consider the bipartite graph $H$ on vertex sets $X=X_1\cup X_2$ and $Y\cup Z$. From its origins in 20 matchings, the maximum degree of $H$ is at most 40, despite the duplication of some vertices.

Now, take any set $Z'\subset Z$ such that $|Z'|=n/3$. We will verify Hall's matching condition between $X$ and $Y\cup Z'$ to show a matching exists, and demonstrate that $H$ satisfies the lemma. For $A\subset X$, let $A'$ be the larger set of $A\cap X_1$ and $A\cap X_2$, so that $|A|\geq|A'|\geq |A|/2$.

Suppose $|A'|\leq m/4$. Then $|N_H(A,Y\cup Z')|\geq|N_H(A',Y)|\geq 2|A'|\geq |A|$, by the properties of $G$.

Suppose $|A'|\geq m/4$ and $|A|\leq n-m/2$. Then there are no edges between $A'$ and $Y\setminus N_H(A')$, so $|Y\setminus N_H(A')|\leq m/4$ and $|Z\setminus N_H(A')|\leq m/4$, using the properties of $G$. Therefore, $|N_H(A',Y\cup Z')|\geq |Y\cup Z'|-m/2\geq |A|$.

Suppose finally then that $|A|\geq n-m/2$. By taking a subset of $A$ of size $n-m/4$ we can see from the previous paragraph that $|N_H(A,Y\cup Z')|\geq n-m/2$. Let $B=(Y\cup Z')\setminus N(A)$, and let $B'$ be the larger set of $B\cap Y$ and $B\cap Z'$, so that $|B'|\leq m/2$. Using similar arguments for $B'$ as we did for $A'$, we can show that $|N_H(B,X)|\geq |N_H(B',X_1)|\geq 2|B'|\geq |B|$. But $N_H(B,X)\subset X\setminus A$, so $|B|\leq|N_H(B,X)|\leq |X\setminus A|=n-|A|$. Therefore $|N_H(A,Y\cup Z')|=n-|B|\geq |A|$, as required.
\oof

\section{Constructing Paths}\label{paths}

In order to build our absorbers in graphs with expansion properties we will need to construct disjoint paths of precise lengths between pairs of vertices. In this section we give such a construction. Similar paths in random graphs were found by Broder, Frieze, Suen and Upfal~\cite{BFSU96}.

\subsection{Between many vertex pairs we can find one path}

\lem \label{connect} Let $m,n\in \N$ satisfy $m\leq n/800$, let $d=n/200m$ and let $n$ be sufficiently large. Let the graph $G$ with $n$ vertices have the property that any set $A\subset V(G)$ with $|A|=m$ satifies $|N(A)|\geq (1-1/64)n$. Suppose $G$ contains disjoint vertex sets $X$, $Y$ and $U$, with $X=\{x_1,\ldots,x_{2m}\}$, $Y=\{y_1,\ldots, y_{2m}\}$ and $|U|=\lceil n/8\rceil$. Suppose, in addition, we have integers $k_i$, $i\in[2m]$, satisfying $2\log n/\log d\leq k_i\leq n/40$. Then, for some $i$, there is an $x_i,y_i$-path of length $k_i$ whose internal vertices lie in $U$.
\ma
\pr

Divide $U$ into two sets, $U_1$ and $U_2$, each of size at least $n/16$. Pick a largest subset $B\subset U_1$ such that $|B|\leq m$ and $|N(B,U_1)|<2d|B|$.  Let $V_1=U_1\setminus B$. Suppose $A\subset V_1$ with $0<|A|\leq m$ has $|N(A,V_1)|<2d|A|$. Then $|N(A\cup B,U_1)|<2d(|A|+|B|)$, so by the choice of $B$ we must have that $|A\cup B|\geq m$. Let $C\subset A\cup B$ have size $|C|=m$. Then $|N(C)|\geq (1-1/64)n$. Since $|A\cup B|\leq 2m$, we have $|N(A\cup B)|\geq (1-1/64)n-m$, so that $|V(G)\setminus (N(A\cup B)\cup A\cup B)|\leq n/64+m$. We thus have,
\[
|N(A,V_1)|\geq |N(A\cup B,U_1)|-|N(B,U_1)|\geq |U_1|-n/64-m-2dm\geq 2dm.
\]
We have, then, that $|V_1|> n/16-m\geq n/20$ and every set $A\subset V_1$ with $|A|\leq m$ satisfies $|N_{V_1}(A)|\geq2d|A|$.

Similarly, find a set $V_2\subset U_2$ with the same expansion property, with $|V_2|\geq n/20$.

Now if $X'\subset X$ is a set of size $m$ then $|N(X',V_1)|\geq |V_1|-n/64\geq n/40$ and hence some vertex $x\in X'$ must have at least $n/40m\geq 2d$ neighbours in the graph $V_1$. Therefore, at least $m+1$ vertices in $X$ have at least $2d$ neighbours in $V_1$.

Similarly at least $m+1$ vertices in $Y$ have at least $2d$ neighbours in $V_2$. Therefore, there is some index $j\in[2m]$ for which $x_j$ and $y_j$ each have at least $2d$ neighbours in $V_1$ and $V_2$ respectively.

The graph $H=G[V_1\cup\{x_j\}]$ then has the property that, given any set $A\subset V(H)$, if $0<|A|\leq m$, then $|N_H(A)|\geq (d+1)|A|+1$, and, if $|A|=m$, then $|N_H(A)|\geq |H|-n/64\geq n/40$. Let $T$ be the $d$-ary tree of depth $l=\lceil \log m/\log d\rceil$. As $k_j\geq 4\log n/\log d$, we have, for sufficiently large $n$, that $\lfloor k_j/2\rfloor-l-1\geq 0$. Attach a path of length $\lfloor k_j/2\rfloor-l-1$ to the root of $T$ to get the tree $T'$ and let the end vertex of the path which is not the root of $T$ be $t_1$. The tree $T'$ has at most $k_j/2+2d^l\leq n/80+2m(d+1)\leq n/40-2md$ vertices. Therefore, by Theorem \ref{almost2}, $H$ contains a copy of $T'$ so that the vertex $t_1$ is embedded onto the vertex $x_j$. Say this copy of $T'$ is $S_1$.

Similarly, $G[V_2\cup\{y_j\}]$ contains a $d$-ary tree with a path of length $\lceil k_j/2 \rceil -l-1$ connecting the root of the regular tree to $y_j$. Call this tree $S_2$.

The set of vertices in the last level of the $d$-regular trees each contain at least $m$ vertices. Suppose these sets are $W_1$ and $W_2$ for the trees $S_1$ and $S_2$ respectively. Let $U'=U\setminus (V(S_1)\cup V(S_2))$, so that $|U'|\geq n/16$. Therefore, $|U'\cap W_1\cap W_2|\geq n/16-2n/64$, so we can pick a vertex $u\in U'$ which is a common neighbour of some vertex $v_1\in W_1$ and some vertex $v_2\in W_2$. Taking the path of length $\lceil k_j/2\rceil-1$ through the tree $S_1$ from $x_j$ to $v_1$, the path $v_1uv_2$ and the path of length $\lfloor k_j/2\rfloor-1$ through the tree $S_2$ from $v_2$ to $y_j$, we get a path from $x_j$ to $y_j$ of length $k_j$, with internal vertices in $U$, as required.
\oof

\subsection{Connecting given vertex pairs}

\lem \label{connectexpand} Let $G$ be a graph with $n$ vertices, where $n$ is sufficiently large, and let $d=160\log^2 n/\log\log n$. Suppose $r,k_1,\ldots,k_r$ are integers with $4\lceil\log n/ \log\log n\rceil\leq k_i\leq n/40$, for each $i$, and $\sum_ik_i\leq 3|W|/4$. Suppose $G$ contains the disjoint vertex pairs $(x_i, y_i)$, $1\leq i\leq r$, and let $W\subset V(G)$ be disjoint from these vertex pairs. 

If $G$ $d$-expands into $W$, then we can find disjoint paths $P_i$, $1\leq i\leq r$, with interior vertices in $W$, so that each path $P_i$ is an $x_i,y_i$-path with length $k_i$.
\ma
\pr Let $k=\lceil\log n/\log\log n\rceil$, and $m=|W|/2d=|W|\log\log n/20\log^2 n$. For any subset $A\subset V(G)$ and $|A|=m$, as $G$ $d$-expands into $W$ and there are no edges between $A$ and $V(G)\setminus (N(A)\cup A)$, we must have that $|N(A)|\geq n-2m\geq (1-1/64)n$.

Using Lemma \ref{splitexpand}, with $n$ sufficiently large, divide $W$ into $k+1$ sets $W_1,\ldots,W_{k}$, and $U$, so that $|W_i|=|W|/16k$, $|U|=15|W|/16$, and $G$ $(2\log n)$-expands into $W_i$, for each $i$. Let $d_0=2\log n$.

We will find as many of the required paths as possible, using vertices from $U$, before finding a $d_0$-matching into $W_1$ from the pairs left to be connected. We can then look at pairs of these neighbours in $W_1$, the connection of any of which will allow us to connect the original pairs. Repeating this with each set $W_i$ in turn will eventually allow us to connect all the pairs of vertices.

Define a \emph{stage-$\a$ situation} to consist of an indexing set $I_\a\subset[r]$, paths $P_i$, $i\in [r]\setminus I_\a$, and sets $T_i,S_i\subset (\cup_{j\leq \a} W_j)\cup\{x_i,y_i\}$ , $i\in I_\a$, as follows.
\begin{itemize}
\item $|I_\a|\leq 2m/(d_0+1)^{\a}$, and $|S_i|=|T_i|=(d_0+1)^{\a}$,
\item The paths, $P_i$, and subsets, $S_i$, $T_i$, are all disjoint,
\item The path $P_i$ is contained in $U\cup(\cup_{j\leq\a}W_j)\cup\{x_i,y_i\}$, has length $k_i$ and end vertices $x_i$ and $y_i$,
\item For each $x\in S_i$, there is an $x,x_i$-path length at most $\a$ in $G[S_i]$, and,
\item For each $y\in T_i$, there is a $y,y_i$-path length at most $\a$ in $G[T_i]$.
\end{itemize}

We will prove the lemma by creating a stage-$k$ situation. Indeed, in such a situation $|I_k|\leq 2m/d_0^k<1$, so that $I_{k}=\emptyset$ and we must have all the required paths.

We can reach a stage-0 situation as follows. Using vertices in $U$ connect as many of the vertex pairs together with paths of the required length, calling these paths $P_i$. Let $I_0$ index the vertex pairs we have yet to connect. Let $U'$ be the vertices in $U$ not covered by any of these paths, $P_i$. Then $|U'|\geq |U|-\sum_ik_i\geq |W|/8$. By Lemma \ref{connect}, applied to the graph $G$ restricted to the set $W$ and the vertex pairs which still need connected, there can be at most $2m$ pairs needing connected, otherwise we could find another path through $U'$ of the correct length. Therefore $|I_0|\leq 2m$. Setting $S_i=\{x_i\}$ and $T_i=\{y_i\}$, for each $i$, gives a stage-0 situation.

We will prove by induction that, given a stage-$\a$ situation with $\a\leq k-1$, we can create a stage-$(\a+1)$ situation. Thus, by induction we can reach a stage-$k$ situation and complete the required paths. 

Suppose then we have a stage-$\a$ situation. We have that $|\cup_{i\in I_\a}(S_i\cup T_i)|=2(d_0+1)^\a|I_\a|\leq 4m\leq |W_{\a+1}|/2d_0$. As $G$ $d_0$-expands into $W_{\a+1}$, we can find a $d_0$-matching from $\cup_{i\in I_\a}(S_i\cup T_i)$ into $W_{\a+1}$. Let $S_i'$ be the union of $S_i$ and its image under this matching. Form $T'_i$ similarly. Then $|T'_i|=|S'_i|=(d_0+1)|S_i|=(d_0+1)^{\a+1}$. The last two properties of a stage-$(\a+1)$ situation are satisfied as $S_{i+1}\subset S_i\cup N(S_i)$ and $T_{i+1}\subset T_i\cup N(T_i)$. 

For each $i\in I_\a$, look for an $x_i$, $y_i$-path of length $k_i$ with interior vertices in $U\cup S_i\cup T_i$ and which is disjoint from all the paths $P_j$ found so far. If such a path exists, take one, and call it $P_i$.

Let $I_{\a+1}\subset I_\a$ index the pairs which still need to be connected after this process. Let $Z=U\setminus(\cup_i V(P_i))$. From the conditions on the integers $k_i$, we have $|Z|\geq |W|/8$. For each $i\in I_{\a+1}$, pair up the vertices in $S_i$ with those in $T_i$. Give each pair, $(s,t)$ say, the integer $k_i-d_{S_i}(s,x_i)-d_{T_i}(t,y_i)\geq k_i-2\a\geq 2k$. If we can find a path between the pairs of this length, then, by taking in addition an $x_i,s$-path length $d_{S_i}(s,x_i)$ in $G[S_i]$ and a $t,y_i$-path length $d_{T_i}(t,y_i)$ in $G[T_i]$, we have an $x_i,y_i$-path, a contradiction. By Lemma \ref{connect} then, we must have at most $2m$ pairs. But we have, in total, $|I_{\a+1}|d_0^{\a+1}$ pairs from the sets $S_i$ and $T_i$, $i\in I_{\a+1}$. Therefore, $|I_{\a+1}|\leq 2m/(d_0+1)^{\a+1}$ as desired.

The indexing set $I_{\a+1}$, paths $P_i$, $i\in [r]\setminus I_{\a+1}$, and sets $S'_i$, $T'_i$, $i\in I_{\a+1}$, form a stage $(\a+1)$-situation, as required.
\oof

\section{Proof of Theorem \ref{mainexpo1}}\label{main}

\de In a graph $G$, we say $(R,r,s)$ is an \emph{absorber} for a vertex $v\in V(G)$ if $\{r,s\}\subset R\subset V(G)$, $r\neq s$, and if there is an $r,s$-path in $G$ with vertex set $R$ and an $r,s$-path in $G$ with vertex set $R\cup\{v\}$. We say that $|R|$ is the \emph{size} of the absorber, and call the vertices $r,s$ the \emph{ends} of the absorber.
\fn

We will use absorbers to develop some flexibility, so that we may finish the embedding of a tree by absorbing unused vertices.

\subsection{Covering Expanders with paths}
We can find absorbers for each vertex in a given set using the following lemma. It uses Lemma \ref{connectexpand}, which says, roughly, that if $G$ $d$-expands into a set $W$, then, for large enough $d$, we can create paths covering up to $3/4$ of the vertices of $W$.

\lem\label{absorbone} Let $n$ be sufficiently large and $d=20\log^2 n$.
Suppose $G$ is a graph containing $W\subset V(G)$, and that $G$ $d$-expands into $W$. 
Then, given any set $A\subset V(G)\setminus W$ with $|A|\leq |W|/300\log^2n$, we can find, disjointly in $G[W]$, 
40 absorbers size $\log^2n+2$ for each vertex $v\in A$.
\ma
\pr Using Lemma \ref{splitexpand}, with $n$ sufficiently large, partition $W$ as $W_1\cup W_2\cup W_3$, so that, for each $i$, $|W_i|=n/3$ and $G$ $(\log^2 n)$-expands into $W_i$. As $|A|\leq |W_1|/160$, we can find an 80-matching from $A$ into $W_1$, by Hall's theorem. Say $v\in A$ is matched to $B_v$ under such a matching, and partition $B_v$ into pairs.

We wish to create $40$ absorbers for each $v\in A$, each using a different pair of vertices from $B_v$. Each absorber will be found by creating paths using $W_2$ and then using $W_3$. As we seek in total to create $40|A|$ absorbers, covering $40|A|(\log^2n+2)\leq 3|W_2|/4=3|W_3|/4$ vertices, we can create all the required absorbers simultaneously and disjointly using Lemma \ref{connectexpand}. To reduce notation, we will describe the construction of one absorber for $v\in A$ using one pair of vertices, $\{x_0,y_1\}\subset B_v$, but, by applying Lemma \ref{connectexpand} to $W_2$ and $W_3$ we can create all the required absorbers at once.

Let $k=\log n$. The graph $G$ $(\log^2n)$-expands into $W_2$. Using Lemma \ref{connectexpand}, we can find a path $Q$ of length $2k+1$ from $x_0$ to $y_1$ with interior vertices in $W_2$. Let $Q$ be the path $x_0x_1x_2\ldots x_{k}y_0y_{k}y_{k-1}\ldots y_2y_1$.

The graph $G[W_3]$ $(\log^2 n)$-expands into $W_3$. By Lemma \ref{connectexpand}, we can find $k$ disjoint paths $P_i$ length $(k-1)$ connecting $x_i$ and $y_i$ pairwise.

\setlength{\unitlength}{0.05cm}
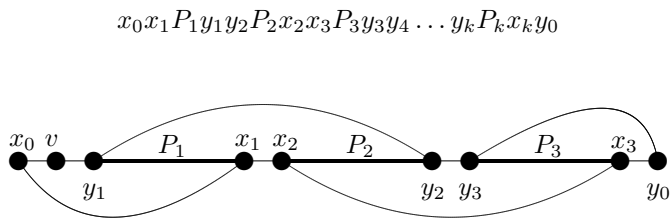
\begin{figure}[b]
\centering
\begin{picture}(210,50)

\put(20,30){\circle*{5}}
\multiput(30,30)(50,0){4}{\circle*{5}}
\multiput(40,30)(50,0){4}{\circle*{5}}

\put(17.5,34){$x_0$}
\put(27,34){$v$}
\put(37,21){$y_1$}
\put(77.5,34){$x_1$}
\put(57,32){$P_1$}
\put(107,32){$P_2$}
\put(157,32){$P_3$}
\put(87.5,34){$x_2$}
\put(177.5,33.3){$x_3$}
\put(127,21){$y_2$}
\put(137,21){$y_3$}
\put(187,21){$y_0$}

\linethickness{0.04cm}
\put(40,30){\line(1,0){40}}
\put(90,30){\line(1,0){40}}
\put(140,30){\line(1,0){40}}

\linethickness{0.01cm}
\qbezier(140, 30)(185, 58)(190, 30)
\qbezier(90, 30)(135, 0)(180, 30)
\qbezier(40, 30)(85, 60)(130, 30)
\qbezier(20, 30)(40, 0)(80, 30)

\put(20,30){\line(1,0){170}}
\end{picture}

\vspace{-0.75cm}
\caption{An absorber for $v$, with much shorter paths depicted.}
\label{absorberpicture}

\end{figure}

Let $R=\cup_iV(P_i)\subset W$. When $k$ is even, the following two $x_0,y_0$-paths have vertex sets $R$ and $R\cup\{v\}$ respectively (see Figure \ref{absorberpicture}, where the heavy lines are paths and the light lines are edges). 
\[
x_0x_1P_1y_1y_2P_2x_2x_3P_3y_3y_4\ldots y_{k}P_{k}x_{k}y_0
\]
\[
x_0vy_1P_1x_1x_2P_2y_2y_3P_3x_3\ldots x_{k}P_{k}y_{k}y_0
\]
When $k$ is odd the following two $x_0,y_0$-paths have vertex sets $R$ and $R\cup\{v\}$ respectively.
\[
x_0x_1P_1y_1y_2P_2x_2x_3P_3y_3y_4\ldots x_{k}P_{k}y_{k}y_0
\]
\[
x_0vy_1P_1x_1x_2P_2y_2y_3P_3x_3\ldots y_{k}P_{k}x_{k}y_0
\]
Hence, $(R,x_0,y_0)$ is an absorber for $v$ in both cases. Counting the vertices created reveals that $|R|=k^2+2$.
\oof

We use our absorbers for single vertices to build up a system of paths with a global absorption property. Suppose we have two absorbers $(S,s,t)$ and $(S',s',t')$, for the vertices $v$ and $v'$ respectively. Given a $t,s'$-path with interior vertices not in $S\cup S'\cup\{v,v'\}$, suppose we add its vertices to $S\cup S'$ to get the set $S''$. It can be easily seen that $(S'',s,t')$ is an absorber both for $v$ and for $v'$. Similarly, given several absorbers with paths linking their endpoints in some order, we can merge them into a single absorber.

The challenge in our setting is that we are concerned with the precise length of the path eventually found through an absorber. We will take our absorbers from Lemma \ref{absorbone} and divide them into groups of up to 40, connecting them as described above. This creates an absorber capable of absorbing up to 40 vertices. When we come to absorb vertices, each path will absorb exactly one vertex so that we know the length of the resulting path. By controlling precisely which absorbers we group together (using Lemma \ref{flexiblematching} for guidance), we can build up some global flexibility.

\lem \label{absorbmany} Let $n$ and $r$ be sufficiently large, and let $l=10^3\log^2n$. Let the graph $G$, with $n$ vertices, contain the disjoint sets $A$, $W$, $X=\{x_1,\ldots,x_{3r}\}$ and $Y=\{y_1,\ldots,y_{3r}\}$. Suppose $r\leq|W|/6l$, $|A|=2r$ and that $G$ $(400\log^2 n)$-expands into $W$. Then, there is a subset $W'\subset W$, with $|W'|=3r(l-2)-r$, as follows. Given any set $A'\subset A$ with $|A'|=r$, there is a set of $3r$ vertex disjoint $x_i,y_i$-paths in $W'\cup A'$ of length $l-1$. In fact, therefore, such paths cover the set $W'\cup A'$.
\ma
\pr
Using Lemma \ref{splitexpand}, with $n$ sufficiently large, partition $W$ as $W_1\cup W_2\cup W_3$, such that, for each $i$, $|W_i|=|W|/3$ and $G$ $(20\log^2 n)$-expands into $W_i$.

Pick a set $B\subset W_1$ with $|B|=|A|$. In the set $W_2$, create the disjoint absorbers $(R_{v,j},r_{v,j},s_{v,j})$, $v\in A\cup B$, $1\leq j\leq 40$, using Lemma \ref{absorbone}, so that each absorber $R_{v,j}$ has size $\log^2n+2$ and can absorb $v$. To describe how to route our paths through these absorbers we refer to an ancillary bipartite graph $H$ with the following properties, provided, for $r$ sufficiently large, by Lemma \ref{flexiblematching}. The bipartite graph $H$ has maximum degree $40$ and vertex classes $[3r]$ and $A\cup B$ such that for any subset $A'\subset A$ with $|A'|=r$ there is a matching between $[3r]$ and $A'\cup B$ in $H$. For each vertex $v\in A\cup B$, let $c_v:N_H(v)\to [40]$ be an injective function, listing the neighbours of $v$ in $H$.

Pick an index $j\in[3r]$, and take the absorbers $R_{v,c_v(j)}$, $v\in N_H(j)$. Take the set of vertices $\{r_{v,c_v(j)},s_{v,c_v(j)}:v\in N(j)\}\cup\{x_j,y_j\}$. Using Lemma \ref{connectexpand} and the set $W_3$, join pairs from these vertices by paths of length at least $\log n$ in such a way as to create an absorber $(S_j,x_j,y_j)$, with $|S_j|=l-1$, which is an absorber for each vertex $v\in N(j)$. Using Lemma \ref{connectexpand}, as $3lr\leq 3|W_3|/4$, we can do this for all indices $j\in [3r]$ simultaneously so that the paths created are disjoint. Note that, because the functions $c_v$ are injective, the absorbers $S_i$ will be disjoint.

Let $W'=(\cup_i(S_i\setminus\{x_i,y_i\})\cup B$. We claim this is such a set as required by the lemma. Indeed, let $A'\subset A$ be any set of size $r$. From the property of the graph $H$, we can find a matching $M$ between $A'\cup B$ and $[3r]$. For each $i\in [3r]$, find the vertex $v\in A'\cup B$ matched to $i$ in $M$ and take find an $x_i,y_i$-path length $l-1$ through $S_i\cup\{v\}$. These paths cover $W'\cup A'$, as required.
\oof

\begin{theorem}\label{expomain2} Let $n$ be sufficiently large and let $l\in \N$ satisfy $l\geq 10^3\log^2n$ and $l|n$. Let a graph $G$ contain $n/l$ disjoint vertex pairs $(x_i,y_i)$ and let $W=V(G)\setminus (\cup_i\{x_i,y_i\})$. Suppose $G$ $d$-expands into $W$, where $d=10^{10}\log^4 n/\log\log n$. Then we can cover $G$ with $n/l$ disjoint paths $P_i$, length $l-1$, so that, for each $i$,  $P_i$ is an $x_i,y_i$-path.
\end{theorem}

\pr[Proof of Theorem \ref{expomain2}]
Let $l_0=10^3\log^2n$, and $d_0=10^8\log^4 n/\log\log n$. We first note that it is sufficient to prove the lemma with $l=l_0$ and $d=d_0$, as follows. Suppose, indeed, the lemma holds for $l_0$ and $d_0$. Take the set $W$ and use Lemma \ref{splitexpand}, with $n$ sufficiently large, to partition $W$ as $W_1\cup W_2\cup W_3$, with $|W_1|=3|W|/4$ and $|W_2|=|W_3|=|W|/8$, such that $G$ $(d/40)$-expands into $W_1$, $W_2$ and $W_3$. Let $l+1=ql_0+r$, where $0\leq r< l_0$ and $q\in \N$. Now, take $n/l$ vertices from $W_2$ and label them $z_1,\ldots,z_{n/l}$. Using Lemma \ref{connectexpand} and the set $W_1$, create disjoint paths of length $\max\{r,\log n\}$ connecting $x_i$ to $z_i$ pairwise. As $\max\{r,\log n\}\leq l/2$, these paths cover at most $n/2\leq 3|W_1|/4$ vertices. If $\log n>r$, then delete $(\log n)-r$ vertices from the end of each of these paths, starting with $z_i$. Let $z_i$, again, be the end vertex of this remaining path for each $i\in[n/l]$. We have now have paths of length $r$ from $x_i$ to $z_i$ for each $i$.

Note that between any two disjoint vertex sets of size $n/d_0$ there must be some edge, since $G$ $d_0$-expands into $W$. Hence there is an edge within any set of $2n/d_0$ vertices. Greedily then, we can take $(n/l)(q-1)$ disjoint edges inside the remaining vertices in $W_1\cup W_2$. Label these edges as $y_{i,j}x_{i,j+1}$, $i\in[n/l]$, $1\leq j< q$, and set $x_{i,1}=z_i$ and $y_{i,q}=y_i$ for each $i$. Let $W'$ be the vertices in $W$ not in any of these edges, or in the paths ending in $z_i$. Then the pairs $\{(x_{i,j},y_{i,j}):i\in[n/l],j\in[q]\}$, the set $W'$ and the graph $G[W'\cup\{x_{i,j},y_{i,j}:i\in[n/l],j\in[q]\}]$ satisfy the conditions of the lemma for $l_0$ and $d_0$. Covering the vertices in $W'$ with paths length $l_0-1$ between these pairs produces the required $x_i,y_i$-paths in the general case.

Suppose then that $l=10^3\log^2n$ and $d=10^8\log^4 n/\log\log n$. Let $m=\lceil n/2d\rceil$, so that any two disjoint vertex sets size $m$ in $G$ have some edge between them. Let $s=n/10^5\log^3 n$.

Using Lemma \ref{splitexpand}, for $n$ sufficiently large, partition $W$ into $W_1$, $W_2$ and $W_3$ as follows. We have $|W_1|=s\log n$, $|W_2|=3s\log n$ and $|W_3|=|W|-2s\log n\geq n/2$, such that $G$ $(d/10)$-expands into $W_3$ and $G$ $(ds\log n/5|W|)$-expands into each set $W_1$ and $W_2$. Note that $ds\log n/5|W|\geq ds\log n/5n\geq 160\log^2 n/\log\log n$ and $s\log n\leq |W_3|/12l$.

By Lemma \ref{absorbmany}, taking $r=2s\log n$ and $A=W_1\cup W_2$, we can find a subset $W'\subset W_3$ of size $(6(l-2)-1)s\log n$ so that, given any subset $Z\subset W_1\cup W_2$ with $|Z|=2s\log n$, there is a set of disjoint $x_i,y_i$-paths, $i\in[6s\log n]$, of length $l-1$ that cover the vertices $W'\cup Z$.

Let $Z_1=W_3\setminus W'$ and let $\a=l/(2\log n+2)$ be an integer\footnote{We assume this for a clear presentation. Where this is not an integer we make a small adjustment to $l$ so that this is true. Such a small adjustment will easily be tolerated by the proof of Lemma \ref{absorbmany}.}. Take distinct vertices $x_{i,j}\in Z_1$, $6s\log n+1\leq i\leq n/l$, $2\leq j\leq \a$, and let $x_{i,1}=x_i$ and $x_{i,\a+1}=y_i$, $6s\log n+1\leq i\leq n/l$. Let $Z_2\subset Z_1$ be the vertices in $Z_1$ not among these labelled vertices. Connect as many pairs $x_{i,j}$ and $x_{i,j+1}$ as possible by disjoint paths of length $2\log n+2$ using vertices from $Z_2$, stopping if no more can be connected or until there are $s$ such vertex pairs remaining. If after this there are $t$ vertex pairs remaining, where $t> s$, then take among them $m$ pairs $(a_i,b_i)$, $i\in [m]$, so that the vertices $a_i$ and $b_i$ are all distinct. Let $Z_3\subset Z_2$ be the set of vertices not covered by any of the paths found so far. Now $|Z_3|=(2\log n+1)t-(|W_1|+|W_2|)/2\geq s\geq 1000m$. Given any set $U\subset V(G)$ with $|U|=m$, $|N(U,Z_3)|\geq |Z_3|-2m\geq (1-1/128)|Z_3|$. Therefore, the graph $G[Z_3\cup(\cup_i\{a_i,b_i\})]$ satisfies the conditions of Lemma \ref{connect}, and so, by that lemma, there is a path of length $2\log n+2$ between one of the pairs $(a_i,b_i)$ in this graph, contradicting $t>s$. The process stops, therefore, with only $s$ pairs remaining.

Say these $s$ unconnected pairs are $(c_i,d_i)$, $i\in[s]$. Let $Z_3$ be the remaining vertices left in $Z_2$ not covered by any paths. Then $|Z_3|=(2\log n+1)s-(|W_1|+|W_2|)/2=s$. Label $Z_3$ as $\{e_{i}:i\in[s]\}$. Recalling that $G$ $(160\log^2 n/\log\log n)$-expands into $W_1$, we can (comfortably) find a generalised matching into $W_1$ from these vertices and the vertices $c_i$ and $d_i$, such that, for each $i\in[s]$, $e_i$ is matched to two vertices, say $e'_i$ and $e''_i$, and $c_i$ and $d_i$ are matched to one vertex each, say $c_i'$ and $d_i'$ respectively.

Use $W_2$ and Lemma \ref{connectexpand} to find disjoint paths length $\log n-1$ connecting $c'_i$ to $e'_i$, and $e''_i$ to $d'_i$. For each $i$, these paths and vertices form paths length $2\log n+2$ from $c_i$ to $d_i$. Thus, we have found paths length $l$ from $x_i$ to $y_i$ passing through $x_{i,2},\ldots,x_{i,\a}$, for each $i>6s\log n$. The paths used $2s\log n=|W_1\cup W_2|/2$ vertices from $W_1\cup W_2$. Let $Z$ be the vertices of $W_1\cup W_2$ not used in any of the paths. Using the properties of the set $W'$, find a set of $6s\log n$ paths, which connect $x_i$ and $y_i$ pairwise for each $i\in[6s\log n]$, have length $l-1$, and cover the set $W'\cup Z$. This completes all the paths and covers the graph.
\oof

\subsection{Proof of Theorem \ref{mainexpo1}}
\pr Let $T$ be a tree with $n$ vertices and maximum degree at most $\Delta$. Take $n$ vertices, and reveal edges independently with probability $\Delta\log^5 n/ 2$ to get the graph $G_1$. By Lemma \ref{randomexpand}, $G_1$ is almost surely an $(n,\Delta\log^4 n/ 20)$-expander. 

Suppose $T$ has at least $n/(4\cdot 10^3\log^2n)\geq n/\log^3 n$ leaves. Remove $n/\log^3 n$ of these leaves from $T$ to get the tree $T'$. By Theorem \ref{almost}, $G_1$ contains a copy of $T'$. Fix this copy of $T'$ and let $M$ be the set of vertices which need leaves attached to them to complete the embedding of $T$. Let the set of vertices not yet in the embedding be $L$. We have $|L|\geq n/\log^3 n$ and $|M|\geq \log^3 n/ \Delta$. Almost surely, by revealing more edges with probability $\Delta\log^5 n/2n$, we can, by Lemma \ref{matchymatchy} find the required generalised matching between $M$ and $L$ to attach the right number of leaves in $L$ to each vertex in $M$ to complete the embedding of $T$.

If $T$ does not have at least $n/(4\cdot 10^3\log^2n)$ leaves, then by Corollary \ref{split}, it must have at least $n/(4\cdot 10^3\log^2n)$ bare paths length $10^3\log^2n$. Remove the interior vertices of such a set of bare paths to get a forest, $T'$, with at most $5n/6$ vertices. Use Lemma \ref{splitexpand}, with $W=V(G)$ and $n$ sufficiently large, to find disjoint sets $W_1$, $W_2$ with $|W_1|=7n/8$ and $|W_2|=n/8$ so that $G$ $(\log^4 n/200)$-expands into $W_i$, $i=1,2$. By Theorem \ref{almost}, $G[W_1]$ contains a copy of $T'$. Fix such a copy and suppose to complete the embedding we need to find paths length $10^3\log^2n$ between $x_i$ and $y_i$ to cover the remaining vertices not in the embedding, $Z$ say. As $W_2\subset Z$, $G[Z\cup(\cup_i\{x_i,y_i\})]$ $d$-expands into $Z$, and so, by Theorem \ref{expomain2}, we can cover the remaining vertices by disjoint paths of the required length to complete the embedding.
\oof



\section{Paths in directed graphs}\label{technical}

We will note here the following variation of Theorem \ref{expomain2}, which will be useful in~\cite{selflove2}. The out-neighbourhood of a set $A$, that is the vertices $y\notin A$ for which there is some vertex $x\in A$ with an edge directed from $x$ to $y$, is denoted $N^+(A)$, and, similarly, the in-neighbourhood of a set $A$ is denoted by $N^-(A)$.
\begin{theorem}\label{pathcoverexpander} Let $n$ be sufficiently large and let $k\in \N$ satisfy $k\geq 10^3\log^3n$ and $k|n$. Let a directed graph $G$ contain $n/k$ disjoint vertex pairs $(x_i,y_i)$ and let $W=V(G)\setminus (\cup_i\{x_i,y_i\})$. Suppose $G$ has the following two properties.
\begin{enumerate}
\item  For any subset $A\subset V(G)$ with $|A|\leq n/2\log^5 n$, $|N^+(A,W)|\geq |A|\log^5 n$, and $|N^-(A,W)|\geq |A|\log^5 n$.

\item Any two disjoint subsets $A,B\subset V(G)$ with $|A|,|B|\geq n/2\log^5 n$ must have a directed edge from $A$ into $B$.
\end{enumerate}
Then we can cover $G$ with $n/k$ paths $P_i$, length $k-1$, so that, for each $i$, $P_i$ is a directed path from $x_i$ to $y_i$.
\end{theorem}

As may be expected, the proof of Theorem \ref{pathcoverexpander} requires a repeat of our methods leading up to the proof of Theorem \ref{expomain2}, except with directed graphs instead of undirected graphs. Typically the expansion conditions are replaced by the same condition on both the out and in neighbourhood of a set. With directed graphs the path connection results we proved (Lemmas \ref{connect} and \ref{connectexpand}) naturally create not just paths of a chosen length, but with chosen orientations too.

For a directed version of Lemma \ref{connect} we require a directed version of Theorem \ref{almost2}. This can be shown to hold, but to reduce the burden of checking, it is possible to use instead a theorem by Friedman and Pippenger~\cite{FP87} whose proof is shorter than that of Theorem \ref{almost2}. Checking through the proof in~\cite{FP87} demonstrates the following holds.
\begin{theorem}\label{fp}
Let $T$ be an oriented tree on $k$ vertices of maximum degree at most $\Delta$, rooted at $r$. Let $H$ be a non-empty directed graph such that, for every $X\subset V(H)$ with $|X|\leq 2k-2$, $|N^+_H(X)|, |N^-_H(X)|\geq(d+1)|X|$ holds. Let $v\in V(H)$ be an arbitrary vertex of $H$. Then $H$ contains a copy of $T$, rooted at $v$.
\end{theorem}

Using Theorem \ref{fp} we can prove versions of Lemmas \ref{connect} and Lemma \ref{connectexpand} in which the orientations of the paths are specified.

The only major change now needed to prove Theorem \ref{pathcoverexpander} is in the construction of an absorber. To achieve a directed absorber, first consider the structure in Figure \ref{reversiblepath}, which we call a \emph{reversible path from $x$ to $y$}. In the diagram, heavy lines denote directed paths and light lines denote edges. Observe that this structure contains directed paths in both directions between $x$ and $y$, both of which use all the vertices. We can now construct a directed absorber as in Figure \ref{absorberpicture2}, where in this case, the heavy lines indicate reversible paths of the kind shown in Figure \ref{reversiblepath}. The absorber contains two directed paths from $x_0$ to $y_0$, one of which contains all the vertices, and the other all the vertices but $v$. The additional construction of reversible paths means that the absorbers will have $\log^3 n$ vertices, instead of $\log^2 n$ vertices. This results in the additional $\log n$ factor in Theorem \ref{pathcoverexpander}. The rest of the proof follows through without any additional complications.

\setlength{\unitlength}{0.05cm}
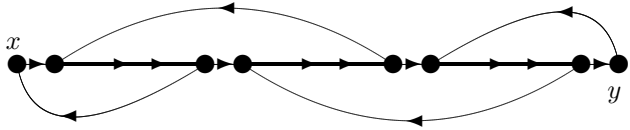
\begin{figure}[h]
\centering
\begin{picture}(210,50)

\multiput(30,30)(50,0){4}{\circle*{5}}
\multiput(40,30)(50,0){4}{\circle*{5}}

\put(27,34){$x$}

\put(187,21){$y$}

\linethickness{0.04cm}
\put(40,30){\line(1,0){40}}
\put(90,30){\line(1,0){40}}
\put(140,30){\line(1,0){40}}

\linethickness{0.01cm}
\qbezier(140, 30)(185, 58)(190, 30)
\qbezier(90, 30)(135, 0)(180, 30)
\qbezier(40, 30)(85, 60)(130, 30)
\qbezier(30, 30)(35, 2)(80, 30)

\put(30,30){\line(1,0){160}}
\linethickness{2cm}
\thicklines
\put(83.5,45){\vector(-1,0){0.000000005}}
\put(133.5,15){\vector(-1,0){0.000000005}}
\put(41.5,16.15){\vector(-1,0){0.000000005}}
\put(173.5,43.95){\vector(-1,0){0.000000005}}

\multiput(37.5,30)(50,0){4}{\vector(1,0){0.000000005}}
\multiput(60,30)(50,0){3}{\vector(1,0){0.000000005}}
\multiput(70,30)(50,0){3}{\vector(1,0){0.000000005}}

\end{picture}

\vspace{-0.75cm}
\caption{A reversible path from $x$ to $y$.}
\label{reversiblepath}

\end{figure}

\setlength{\unitlength}{0.05cm}
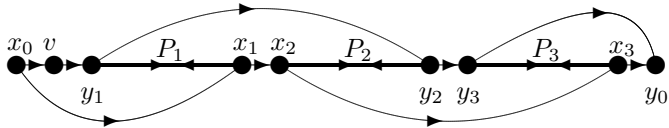
\begin{figure}[h]
\centering
\begin{picture}(210,50)

\put(20,30){\circle*{5}}
\multiput(30,30)(50,0){4}{\circle*{5}}
\multiput(40,30)(50,0){4}{\circle*{5}}

\put(17.5,34){$x_0$}
\put(27,34){$v$}
\put(37,21){$y_1$}
\put(77.5,34){$x_1$}
\put(57,32){$P_1$}
\put(107,32){$P_2$}
\put(157,32){$P_3$}
\put(87.5,34){$x_2$}
\put(177.5,33.3){$x_3$}
\put(127,21){$y_2$}
\put(137,21){$y_3$}
\put(187,21){$y_0$}

\linethickness{0.04cm}
\put(40,30){\line(1,0){40}}
\put(90,30){\line(1,0){40}}
\put(140,30){\line(1,0){40}}

\linethickness{0.01cm}
\qbezier(140, 30)(185, 58)(190, 30)
\qbezier(90, 30)(135, 0)(180, 30)
\qbezier(40, 30)(85, 60)(130, 30)
\qbezier(20, 30)(40, 0)(80, 30)

\put(20,30){\line(1,0){170}}
\linethickness{2cm}
\thicklines
\put(83.5,45){\vector(1,0){0.000000005}}
\put(133.5,15){\vector(1,0){0.000000005}}
\put(47.5,15.15){\vector(1,0){0.000000005}}
\put(177.5,43.95){\vector(1,0){0.000000005}}

\put(27.5,30){\vector(1,0){0.000000005}}
\multiput(37.5,30)(50,0){4}{\vector(1,0){0.000000005}}
\multiput(60,30)(50,0){3}{\vector(1,0){0.000000005}}
\multiput(63,30)(50,0){3}{\vector(-1,0){0.000000005}}
\end{picture}

\vspace{-0.75cm}
\caption{An absorber for $v$, for directed graphs. The paths $P_i$ are reversible.}
\label{absorberpicture2}

\end{figure}

\section{Future work}\label{future}

In the proof of Theorem \ref{mainexpo1}, we revealed the graph entirely before we embedded trees with long bare paths. Therefore, this proof gives a ``universal'' result for trees of this kind, that is, almost surely the random graph contains all such trees simultaneously. For trees with many leaves we revealed more edges in the graph after we had embedded part of the tree, so we don't have a universal result. Johannsen, Krivelevich and Samotij~\cite{JKS12} showed that there is some constant $c>0$ for which, if $\Delta\geq \log n$, the random graph $\GG(n,c\Delta n^{-1/3}\log n)$ almost surely contains every tree in $\TT(n,\Delta)$. Ferber, Nenadov and Peter~\cite{FNP13} recently showed that the random graph $\GG(n,p)$ almost surely contains every tree in $\TT(n,\Delta)$ if $p=\omega(\Delta^8 n^{-1/2}\log^2 n)$. By giving a construction which can embed all trees with many leaves, and by improving the methods used here slightly, we will show that for every $\Delta>0$ there is some $f(\Delta)$ so that the random graph $\GG(n,f(\Delta)\log^2 n /n)$ is almost surely $\TT(n,\Delta)$-universal~\cite{selflove3}. 

The probability required in~\cite{selflove3} seems to be an artifact of the construction of absorbers given here. It seems plausible that the threshold for a random graph to contain all the trees in $\TT(n,\Delta)$, for fixed $\Delta$, is the same as that conjectured by Kahn for single spanning trees (see Section \ref{intro}). As mentioned in Section \ref{intro}, the author anticipates giving a proof of the conjecture for single trees. However, the techniques will not apply to the more difficult problem of constructing all the trees in $\TT(n,\Delta)$ simultaneously.

\begin{acknowledgements} The author would like to thank Andrew Thomason for his help and suggestions.
\end{acknowledgements}

\bibliographystyle{plain}
\bibliography{rhmreferences}

\end{document}